\documentclass[twoside]{article}
\usepackage{amsmath}
\usepackage{amssymb}
\usepackage{amsthm}
\usepackage{latexsym}
\usepackage{hyperref}

\newcommand{\mb}{\mathbb}
\newcommand{\mr}{\mathrm}
\newcommand{\ma}{\mathcal}
\newcommand{\mf}{\mathfrak}
\newcommand{\ts}{\textstyle}
\newcommand{\hot}{\widehat{\otimes}}
\newcommand{\inp}[2]{\langle #1 \,,\, #2 \rangle}

\newtheorem{thm}{Theorem}
\newtheorem{cor}[thm]{Corollary}
\newtheorem{lem}[thm]{Lemma}
\newtheorem{prop}[thm]{Proposition}
\newtheorem{defn}[thm]{Definition}

\begin{document}
\begin{center}
\huge Some Fr\'echet algebras for which the\\ 
Chern character is an isomorphism\\[5mm]
\Large Maarten Solleveld\\
\large October 2005
\\[1cm]
\end{center}

\begin{minipage}{13cm}
\textbf{Abstract.}
Using similarities between topological $K$-theory and periodic cyclic 
homology we show that, after tensoring with $\mb C$, for certain 
Fr\'echet algebras the Chern character provides an isomorphism 
between these functors. This is applied to prove that the Hecke algebra 
and the Schwartz algebra of a reductive $p$-adic group have isomorphic 
periodic cyclic homology.\\
Later an appendix was added, to deal with infinite direct products of algebras.  
\end{minipage}
\\[1cm]

\noindent
Let $X$ be a smooth manifold. The classical Chern character is a 
map that assigns to a vector bundle on $X$ a class in the even De Rham
cohomology of $X$. This extends naturally to a ring homomorphism
\begin{equation}
\mr{Ch} : K^*(X) \to H_{DR}^*(X) 
\end{equation}
More or less since the beginning of topological $K$-theory \cite{AH}  
it has been known that if $X$ is compact this yields an 
isomorphism of $\mb Z / 2 \mb Z$-graded algebras
\begin{equation}\label{eq:1}
\mr{Ch} \otimes \mr{id} : K^*(X) \otimes \mb R \to H_{DR}^*(X)
\end{equation}

However for noncompact $X$ the graded vector spaces in \eqref{eq:1}
are not necessarily isomorphic. This is because the $K$-theory of
$X$ is not defined directly: one first has to take the one-point
compactification of $X$, then determine the $K$-groups of that space, 
and finally take the quotient by a suitable subgroup isomorphic to
$K_*$(point). In other words topological 
$K$-theory, contrarily to De Rham cohomology, is similar to singular 
(or \u Cech) cohomology \emph{with compact support.}
 
If we extend our coefficients from $\mb R$ to $\mb C$ then both sides
of \eqref{eq:1} can be expressed in terms of the Fr\'echet-algebra
$C^\infty (X)$ of complex-valued smooth functions on $X$:
\begin{align}
\label{eq:23} K^*(X) & \cong K_*(C^\infty (X))\\
\label{eq:15} H_{DR}^*(X ; \mb C) & \cong H^{DR}_* (C^\infty (X))
\end{align}
Here \eqref{eq:23} follows from the smooth version of the Serre-Swan
theorem, while \eqref{eq:15} is a consequence of the canonical 
identification
\begin{equation}\label{eq:22}
\Omega^*(X;\mb C) = \Omega_* (C^\infty(X))
\end{equation}
of the De Rham complex over $X$ with the exterior algebra of 
$C^\infty(X)$.

From the above we get an isomorphism
\begin{equation}\label{eq:2}
\mr{Ch} \otimes \mr{id} : K_*(A) \otimes \mb C \to H^{DR}_*(A)
\end{equation} 
for $A = C^\infty (X)$. Clearly it is desirable to extend 
\eqref{eq:2} to other (noncommutative) Fr\'echet algebras $A$. 
Indeed this is a well-studied subject in noncommutative geometry, 
and strong results have been achieved. The first complication on 
this path is that, although it is always defined, De Rham 
(co)homology behaves well mainly for commutative algebras. Therefore 
we must put another functor on the right hand side of \eqref{eq:2}, 
and the best choice turns out to be the topological periodic cyclic 
homology $HP_*$. This was justified by Connes, who showed \cite{Con} 
that in the above case we have
\begin{equation}\label{eq:3}
HP_*(C^\infty (X)) \cong H^{DR}_* (C^\infty (X))
\end{equation}
Cyclic homology is intimately related to Hochschild homology, analogous
to the relation between \eqref{eq:15} and \eqref{eq:22}.
Indeed to derive \eqref{eq:3} Connes first proved a smooth version of 
the (algebraic) Hochschild-Kostant-Rosenberg theorem :
\begin{equation}\label{eq:24}
HH_*(C^\infty(X)) \cong \Omega^*(X;\mb C)
\end{equation}

Before we continue we have to specify in exactly which category of 
algebras we wish to work, because this matters for the definitions and 
properties of the functors $K_*$ and $HP_*$. We use the largest that 
seems reasonable (at present).

Recall that an $m$-algebra is a complete locally convex algebra whose
topology is defined by a family of submultiplicative seminorms. 
It can be shown that $m$-algebras are precisely the projective limits 
of Banach algebras, and in particular every $C^*$-algebra is an 
$m$-algebra. Fr\'echet $m$-algebras can be described as those 
$m$-algebras for which a countable number of seminorms suffices to 
define the topology, or equivalently as those Fr\'echet algebras whose
seminorms may be taken submultiplicative. Unless explicitly stated 
otherwise, by the (topological) tensor product of two $m$-algebras we 
will always mean  the projective tensor product, which is also an 
$m$-algebra.   

The $K$-theory of $m$-algebras was defined by Phillips \cite{Phi} 
(for Fr\'echet $m$-algebras) and by Cuntz \cite{Cun2}, as a 
special case of a bivariant functor on $m$-algebras. 
It satisfies the following properties:
\begin{enumerate}
\item additivity: $K_i \big( \prod_{n=1}^\infty A_n \big) 
      \cong \prod_{n=1}^\infty K_i(A_n)$
\item stability: $K_i(M_n(A)) \cong K_i(A)$
\item 2-periodicity: $K_{i+2}(A) \cong K_i(A)$
\item diffeotopy-invariance: if $(\phi_t)_{t \in [0,1]}$ is a diffeotopy
      of $m$-algebra homomorphisms then $K_i(\phi_0) = K_i(\phi_1)$
\item excision: if $I$ is an ideal of $A$ then there exists an exact hexagon
\begin{center}
$\begin{array}{ccccc}
K_0(I) & \to & K_0(A) & \to & K_0(A/I) \\
\uparrow & & & & \downarrow \\
K_1(A/I) & \leftarrow & K_1(A) & \leftarrow & K_1(I)
\end{array}$
\end{center}
\item continuity: $K_i \big(\lim \limits_{n \to \infty} 
      A_n \big) \cong \lim \limits_{n \to \infty} K_i(A_n)$, at least for Banach algebras
\end{enumerate}
Of course this list is by no means exhaustive, for example stability as
above is just a weak version of the Morita-invariance of $K$-theory.

Periodic cyclic (co-)homology was first defined by Connes \cite{Con}
for general locally convex algebras, and he showed directly that it 
is invariant under diffeotopies. This homology theory satisfies 3 by 
construction, while 1 and 2 follow from the well known corresponding 
properties of Hochschild homology. Cuntz \cite{Cun2} 
proved that there exists a unique functorial Chern character on
the category of $m$-algebras
\[ch : K_* \to HP_*
\]
which is compatible with 1--4. 

However there are some problems with excision for $HP_*$, and they 
stem from the fact that a closed subspace of a topological vector 
space does not always possess a closed complement. Therefore we 
restrict ourselves to admissible extensions of $m$-algebras, 
i.e. those short exact sequences $0 \to A \to B \to C \to 0$
which admit a continuous linear splitting $C \to B$. Equivalently we 
can require that there exists a continuous projection of $B$ onto $A$.
In the same spirit we call an ideal $I$ of $A$ admissible if the 
associated sequence $0 \to I \to A \to A/I \to 0$ is. 

\begin{thm}\label{thm:1}
Let $0 \to A \to B \to C \to 0$ be an admissible extension of 
$m$-algebras. Then excision holds for $HP_*$ and the various Chern
characters make a commutative diagram
\begin{center}
$\begin{array}{ccccccccccc}
K_1(A) & \to & K_1(B) & \to & K_1(C) & \to & K_0(A) & \to & 
 K_0(B) & \to & K_0(C) \\
\downarrow & & \downarrow & & \downarrow & & \downarrow & & \downarrow
 & & \downarrow \\
HP_1(A) & \to & HP_1(B) & \to & HP_1(C) & \to & HP_0(A) & \to &  
 HP_0(B) & \to & HP_0(C) 
\end{array}$
\end{center}
Moreover if $\eta: K_0(C) \to K_1(A)$ and 
$\partial: HP_0(C) \to HP_1(A)$ are the connecting maps, then\\
ch $\circ \, \eta = 2 \pi i \partial \, \circ$ ch.
\end{thm}
\emph{Proof.}
The excision part is due to Cuntz \cite{Cun1}, while Nistor 
\cite{Nis1} showed that the Chern character commutes with the index 
maps and is compatible with the connecting maps, in the specified 
sense. $\qquad \Box$
\\[2mm] 

Neither is $HP_*$ continuous in general, because 
inductive limits and projective tensor products do not commute. 
On the other hand, inductive limits \emph{do} commute with inductive 
tensor products, so the periodic cyclic homology $HP'_*$ based on 
this tensor product has a better chance of being continuous. 
Indeed Brodzki and Plymen \cite{BP1} showed that if there exists an 
$N \in \mb N$ such that $HH'_i(A_n) = 0 \; \forall i > N$, then 
\begin{equation}\label{eq:20}
HP'_i \big( \lim_{n \to \infty} A_n \big) \cong \lim_{n \to \infty} HP'_i(A_n)
\end{equation}
However, as long as we are not working with nuclear Fr\'echet algebras,
for which there is only one topological tensor product, excision 
remains a bit of a problem for this cyclic theory.

Now let us define the category $\ma C$ of algebras that we are 
going to study :

\begin{defn}\label{def:2}
The category $\ma C$ is a full subcategory of the category of 
$m$-algebras. Its objects are those $m$-algebras $A$ for which
the Chern character induces an isomorphism
\[ch \otimes \mr{id} : K_*(A) \otimes \mb C \to HP_*(A) \otimes_{\mb C} \mb C = HP_* (A)
\]
\end{defn}

\begin{cor}\label{cor:3}
The category $\ma C$ is stable, closed under diffeotopy-equivalences
and under finite direct sums. It contains all algebras of the type 
$C^\infty (X)$, where $X$ is a compact smooth manifold. 

If $0 \to A \to B \to C \to 0$ is an admissible 
extension of $m$-algebras and two elements of $\{A,B,C\}$ are 
objects of $\ma C$, then so is the third.
\end{cor}
\emph{Proof.}
We only prove the last statement, everything else following directly 
from the above remarks. Tensor the commutative diagram of Theorem 
\ref{thm:1} with $\mb C$. If for example (the other cases being 
similar) we know that $\{A,C\} \subset \mr{Ob}(\ma C)$, then first 
consider the diagram obtained by deleting the column containing 
$HP_1(B)$. The five lemma shows that 
$K_0(B) \otimes \mb C \cong HP_0(B)$. Likewise, if we 
retain this column but delete the one with $HP_0(B)$ then we deduce 
that also $K_1(B) \otimes \mb C \cong HP_1(B)$, so that indeed 
$B \in \mr{Ob}(\ma C). \qquad \Box$ 
\\[2mm]

This last property can also be formulated by saying that $\ma C$ is 
closed under taking extensions, ideals and quotients. This can
easily be extended to longer sequences of $m$-algebras.
 
\begin{lem}\label{lem:4}
Consider an increasing sequence
\[0 = I_0 \subset I_1 \subset \cdots \subset I_n \subset I_{n+1} = A
\]
of admissible ideals of $A$. If all the quotients $I_m / I_{m-1}$ 
are objects of $\ma C$, then so are all the $I_m$ themselves.
Similarly let
\[0 \to B_1 \to B_2 \to \cdots \to B_n \to 0
\]
be an admissible exact sequence of $m$-algebras.
If all but one of the $B_i$ belong to $\ma C$, then so does the last.
\end{lem}
\emph{Proof.}
Consider the admissible extensions
\[\begin{array}{ccccccccc}
0 & \to & I_{m-1} & \to & I_m & \to & I_m / I_{m-1} & \to & 0 \\
0 & \to & \mr{im}\: (B_{m-1} \to B_m) & \to & B_m & \to & 
  \mr{im}\: (B_m \to B_{m+1}) & \to & 0  
\end{array}
\]
They degenerate for $m=1$, so the lemma follows from Corollary 
\ref{cor:3}, with induction to $n. \qquad \Box$
\\[2mm]

From now on let $G$ be a finite group acting by diffeomorphisms on 
$X$, which is assumed second countable but not necessarily compact. 
Wassermann \cite{Was2} extended Connes' result \eqref{eq:3} to 
this equivariant noncompact setting by showing that
\begin{equation}\label{eq:4}
HP_*(C^\infty (X)^G) \cong HP_*(C^\infty (X))^G 
\cong H^{DR}_*(C^\infty (X))^G
\end{equation}
By \eqref{eq:15} the right hand side can be identified with the 
$G$-invariant part $H_{DR}^*(X; \mb C)^G$ of the De Rham cohomology 
of $X$, and Grothendieck proved in Corollaire 5.2.3 of 
\cite{Gro} that this in turn is naturally isomorphic to 
$H^*(X/G ; \underline{\mb C})$, the cohomology of the constant sheaf
$\underline{\mb C}$ over the orbifold $X/G$. Since $C^\infty (X)^G$ 
is by definition the algebra of smooth functions on $X/G$, 
\eqref{eq:4} can be restated by saying that \eqref{eq:3} also holds 
for orbifolds that are quotients of manifolds by finite groups.

A related interesting algebra is the crossed product
$C^\infty(X) \rtimes G$. To formulate the relevant results we 
introduce
\[\hat X := \{ (g,x) \in G \times X : g x = x \}
\]
with the $G$-action 
\[g \cdot (g',x) = (g g' g^{-1}, g x)
\] 
The space $\hat X/G$ is called the extended quotient of $X$ by $G$. 
It is a disjoint union of orbifolds, one for each conjugacy class in 
$G$, and (from a homological point of view) well-suited to deal with 
singularities of the group action.
Brylinski \cite{Bry} proved that 
\begin{align}
\label{eq:16} HP_*(C^\infty(X) \rtimes G) &\cong 
\big( H^*_{DR}(\hat X;\mb C) \big)^G \\
\label{eq:25} HH_*(C^\infty(X) \rtimes G) &\cong 
\big( \Omega^*(\hat X;\mb C) \big)^G
\end{align}
By the above result of Grothendieck, \eqref{eq:16} is isomorphic to 
$H^*(\hat X /G; \underline{\mb C})$. Baum and Connes \cite{BC} 
constructed (for compact $X$) an equivariant Chern character 
\begin{equation}\label{eq:17}
\mr{Ch}_G : K_*(C^\infty(X) \rtimes G) \cong K^*_G(X) \to  
\big( H^*_{DR}(\hat X;\mb C) \big)^G
\end{equation}
and proved that it becomes an isomorphism after tensoring with 
$\mb C$. So in this case it has not only been known for long that 
$K_* \otimes \mb C$ and $HP_*$ agree, they have also been determined
in geometrical terms.

By the way, the isomorphisms \eqref{eq:4} till \eqref{eq:17} 
also hold for the algebra $C_c^\infty(X)$ of compactly
supported smooth functions on $X$, provided that one takes De Rham 
cohomology with compact support everywhere. See \cite{BN} for these
and more general results on the cyclic homology of algebras associated 
to orbifolds. 

But things cannot always be this nice. 
In \cite{Was2} it was noticed that, even for compact $X$,
\[HH_*(C^\infty (X)^G) \quad \mr{and} \quad \Omega^*(X;\mb C)^G  
\]
are not isomorphic in general.

If $\mb C^N = \mb C G$ is the regular representation of $G$, then
we can endow $C^\infty( X; \mr{End}(\mb C G)) = 
C^\infty(X) \otimes \mr{End}(\mb C G)$ with the diagonal $G$-action,
and it is well-known that
\begin{equation}
C^\infty(X) \rtimes G \cong C^\infty( X; \mr{End}(\mb C G))^G =
M_N (C^\infty(X))^G
\end{equation}
In our main theorem will show that algebras of this type (and 
some others as well) belong to the category $\ma C$ defined above. 
Note that this does not contradict the remark after equation 
\eqref{eq:1}, since over there we were actually dealing with the
$C^*$-algebra $C_0(X)$, which is quite different from 
$C^\infty (X)$. 

To include manifolds with boundary we recall 
the following conventions :

\begin{defn}\label{def:5}
Let $Y \subset X$ be arbitrary subsets of a smooth manifold $M$, and let $V$ 
be a complex vector space.
\begin{align*}
C^\infty (X) & := \{ f \big|_X : f \in C^\infty(U) \:
  \mbox{for some open U} \,,\: X \subset U \subset M \} \\
C_0^\infty (X,Y) & := \{f \in C^\infty (X) : f \big|_Y = 0 \} \\
C_0^\infty (X,Y ; V) & := C_0^\infty (X,Y) \otimes V
\end{align*}
\end{defn}

\begin{thm}\label{thm:6}
Let $X$ be a smooth manifold with boundary, $N \in \mb N$ and consider 
the $m$-algebra $A = C^\infty (X ; M_N(\mb C)) = M_N(C^\infty(X))$. 
Suppose that a finite group $G$ acts on $A$ by
\begin{equation}\label{eq:6}
g \cdot a(x) = u_g(x) a(\alpha_g^{-1} x) u_g^{-1}(x)
\end{equation}
where $\alpha_g$ is a diffeomorphism of $X$ and $u_g \in A$. Then
$A^G \in \mr{Ob}(\ma C)$, and $K_*(A^G)$ is a finitely generated 
abelian group whenever $X$ is compact.
\end{thm}

First we prove a special case, an equivariant version of the 
Poincar\'e lemma :

\begin{lem}\label{lem:7}
In the setting of Theorem \ref{thm:6}, suppose that $X$ is 
$G$-equivariantly contractible to a point $x_0 \in X$. Then $A^G$ 
is diffeotopy-equivalent to its fiber $\mr{End}_G (\mb C^N)$ over 
$x_0$. In particular $K_*(A^G) = K_0(A^G)$ is a free abelian group of 
finite rank and $A^G \in \mr{Ob}(\ma C)$.
\end{lem}
\emph{Proof.}
Our main task is to adjust the $u_g$ suitably. Since $X$ is 
contractible we can find for every $g \in G$ a function 
$f_g \in C^\infty(X)$ such that $f_g^{-N} = \det(u_g)$. The $G$-action
does not change when we replace $u_g$ with $f_g u_g$, so we may assume
that $\det(u_g) \equiv 1 ,\, \forall g \in G$. The premise that
\eqref{eq:6} is a group action guarantees that there is a smooth 
function $\lambda : G \times G \times X \to \mb C$ such that
\begin{equation}\label{eq:7}
u_g(x) u_h( \alpha_g^{-1} x) = \lambda (g,h,x) u_{gh}(x)
\end{equation}
Taking determinants we see that in fact $\lambda (g,h,x)^N \equiv 1$, 
so $\lambda$ does not depend on $x \in X$. All the elements of 
$\alpha (G)$ fix $x_0$, so the fiber $V_0 = \mb C^N$ over that point 
is a projective $G$-representation $(\pi_0, V_0)$. Thus we are in a 
position to apply Schur's theorem \cite{Sch}, which says that there 
exists a finite central extension $G^*$ of $G$ such that $\pi_0$ lifts 
to a representation of $G^*$. This lift only involves scalar multiples 
of the $u_g(x_0)$, so it immediately extends to $X$. Then 
\eqref{eq:7} becomes the cocycle relation
\begin{equation}\label{eq:8}
u_{gh}(x) = u_g(x) u_h( \alpha_g^{-1} x)
\end{equation}
Notice that still $A^{G^*} = A^G$, so without loss of generality we can
replace $G$ by $G^*$.

Now we want to make the $u_g(x)$ independent of $x \in X$. Wassermann 
\cite{Was1} indicated how this can be done in the continuous case,
and his argument can easily be adapted to our smooth setting.
The crucial observation, first made by Rosenberg \cite{Ros}, is that 
$A^G$ can be rewritten as the image of an idempotent in a larger algebra. 
This idempotent can then be deformed to one independent of $x$.

Indeed, let $A \rtimes_\alpha G$ be the crossed product of $A$ and 
$G$ with respect to the action $\alpha$ of $G$ on $X$, and 
$(r_t)_{t \in [0,1]}$ a smooth $G$-equivariant contraction from $X$ to 
$x_0$. (For smooth manifolds the existence of a continuous contraction 
implies the existence of a smooth one.) Define
\begin{equation}
p_t(x) := |G|^{-1} {\ts \sum_{g \in G}} u_g (r_t x) g
\end{equation}
Then $p_t \in A \rtimes_\alpha G$ is an idempotent by \eqref{eq:8}, 
and by \cite{Ros}
\begin{equation}
\phi_1 : A^G \to p_1 (A \rtimes_\alpha G) p_1 : \sigma \to p_1 \sigma p_1
\end{equation}
is an isomorphism of Fr\'echet algebras. Clearly the idempotents $p_t$ 
are all homotopic, so they are conjugate in the Banach completion 
$C(X;M_N(\mb C)) \rtimes_\alpha G$ of $A \rtimes_\alpha G$. Moreover the 
standard argument for this, as for example in Proposition 4.3.2 of 
\cite{Bla}, shows that $p_0$ and $p_1$ are conjugate by an element of 
$A \rtimes_\alpha G$. Alternatively we can use the stronger result that 
homotopic idempotents in Fr\'echet $m$-algebras are conjugate, but this 
statement is vastly more difficult to prove than its Banach algebra 
version, cf. \cite{Phi}, Lemma 1.12 and Lemma 1.15. In any case, we have 
\begin{equation}\label{eq:9}
A^G \cong p_1 (A \rtimes_\alpha G) p_1 \cong p_0 (A \rtimes_\alpha G) p_0 
\cong C^\infty (X ; \mr{End}_{\mb C} (V_0))^G
\end{equation}
To this last algebra we can apply the obvious diffeotopy 
$\sigma \to \sigma \circ r_t$. This shows that $A^G$ is 
diffeotopy-equivalent to its fiber $\mr{End}_G (V_0)$ over $x_0$, 
and the remaining statements on $K_*(A^G)$ and $HP_*(A^G)$ follow 
from the semisimplicity of the finite dimensional algebra 
$\mr{End}_G (V_0). \quad \Box$
\\[2mm]

To make full use of the proof of this lemma we have to study certain
ideals in such algebras as well :

\begin{lem}\label{lem:8}
Let $U \subset \mb R^n$ be an open bounded star-shaped set. 
The $m$-algebra $C_0^\infty (\mb R^n, \mb R^n \setminus U)$ 
belongs to $\ma C$ and 
\[K_*(C_0^\infty (\mb R^n, \mb R^n \setminus U)) \cong 
  H^*_c (\mb R^n) \cong \mb Z
\] 
is concentrated in degree $n$.
\end{lem}
\emph{Proof.}
Clearly we may assume that 0 is the center of $U$. Let $P$ be 
the point of the $n$-sphere corresponding to infinity under 
the stereographic projection $S^n \to \mb R^n$. By assumption 
$C_0^\infty (\mb R^n, \mb R^n \setminus U) \cong C_0^\infty (S^n, Y)$ 
for some closed neighborhood $Y$ of $P$, 
and we will show that the latter algebra is diffeotopy-equivalent to 
$C_0^\infty (S^n, P)$. Let $(r_t)_{t \in [0,1]}$ be a diffeotopy of
smooth maps $S^n \to S^n$ such that
\begin{enumerate}
\item $\forall t \;: r_t(P) = P$ and $r_t (Y) \subset Y$ 
\item a neighborhood of $-P$ is fixed pointwise by all $r_t$
\item $r_1 = \mr{id}_{S^n}$ and $r_0 (Y) = P$
\end{enumerate}
To construct such maps, we can require that $r_t$ stabilizes every 
geodesic from $-P$ to $P$ and declare that furthermore $r_t(Q)$ 
depends only on $t$ and on the distance from $Q$ to $P$. Then we 
only have to pick a suitable smooth function of $t$ and this 
distance. Given this, consider the $m$-algebra homomorphisms 
\begin{eqnarray*}
\phi : C_0^\infty (S^n, P) & \to & C_0^\infty (S^n, Y) :
  f \to f \circ r_0\\
i : C_0^\infty (S^n, Y) & \to & C_0^\infty (S^n, P) : f \to f
\end{eqnarray*}
By construction $\phi \circ i$ and $i \circ \phi$ are diffeotopic
to the respective identity maps on $C_0^\infty (S^n, Y)$ and  
$C_0^\infty (S^n, P)$, so these algebras are indeed 
diffeotopy-equivalent. Thus we reduced our task to calculating the 
$K$-groups and periodic cyclic homology of $C_0^\infty (S^n, P)$. 
Fortunately there is an obvious split extension
\[0 \to C_0^\infty (S^n, P) \to C^\infty (S^n) \to \mb C \to 0
\]
which by Corollary \ref{cor:3} consists entirely of $m$-algebras
in the category $\ma C$. It is well known that
\[K_*(C^\infty(S^n)) \cong K^*(S^n) \cong \mb Z^2
\]
with one copy of $\mb Z$ in degree 0 and the other in degree $n$.
Since $K_*(\mb C) = K_0(\mb C) \cong \mb Z$ the lemma follows 
from the excision property of the $K$-functor. $\quad \Box$
\\[2mm]

\noindent\emph{Proof of Theorem \ref{thm:6}.}
All our arguments will depend on the existence of a specific cover 
of $X$. To construct it we use a theorem of Illman \cite{Ill}, 
which states that $X$ has a smooth equivariant triangulation. 
In slightly more down-to-earth language this means (among others) 
that there exists a countable, locally finite simplicial complex 
$\Sigma$ in a finite dimensional orthogonal representation space $V$ 
of $G$, and a $G$-equivariant homeomorphism $\psi: \Sigma \to X$. 
Moreover $\psi$ is smooth as a map from a subset of $V$ to $X$, and 
its restriction to any simplex $\sigma$ of $\Sigma$ is an embedding.

For every such $\sigma$ we put
\begin{equation}\label{eq:18}
U'_\sigma := \{ v \in \Sigma : d(v, \sigma) \leq r_\sigma \}
\end{equation}
where $d$ is the Euclidean distance in $V$. We require that the
radius $r_\sigma$ depends only on the $G$-orbit of $\sigma$ and that
$r_\tau > r_\sigma > 0$ if $\tau$ is a face of $\sigma$. The 
orthogonality of the action of $G$ on $V$ guarantees that 
\[g U'_\sigma = U'_{g \sigma} \quad \mr{and} \quad 
U'_\sigma \cap U'_\tau \subset U'_{\sigma \cap \tau}
\]
if we take our radii small enough. Let $D'_\sigma$ be the union, 
over all faces $\tau$ of $\sigma$, of the $U'_\tau$, and $G_\sigma$ 
the stabilizer of $\sigma$ in $G$. From the above we deduce that 
$U'_\sigma \setminus D'_\sigma$ is $G_\sigma$-equivariantly 
retractible to $\sigma \setminus D'_\sigma$.

Now we abbreviate $U_\sigma := \psi (U'_\sigma)$ and 
$D_\sigma := \psi (D'_\sigma)$, so that
$\ma U := \{ U_\sigma : \sigma$ simplex of $\Sigma \}$
is a closed $G$-equivariant cover of $X$.
Let $X_m$ be the union of all those $U_\sigma$ for which 
$m + \dim \sigma \leq \dim X$. It is a closed subvariety (with 
boundary and corners) of $X$ and it is stable under the action 
of $G$. Define the following $G$-stable ideals of $A$ :
\[I_m := \{ a \in A : a \big|_{X_m} = 0 \}
\]
By Th\'eor\`eme IX.4.3 of \cite{Tou}
\begin{equation}
0 \to I_m \to A = C^\infty (X ; M_N( \mb C)) \to 
  C^\infty(X_m ; M_N(\mb C)) \to 0
\end{equation}
is an admissible extension of Fr\'echet algebras. 
Using the finiteness of $G$ we see that $I_m^G$ is 
an admissible ideal in $I_{m+1}^G$ and that
\begin{equation}\label{eq:10}
I_{m+1}^G / I_m^G \cong (I_{m+1} / I_m)^G \cong 
  C_0^\infty (X_m , X_{m+1} ; M_N (\mb C))^G 
\end{equation}
In order to apply Lemma \ref{lem:4} to the sequence
\begin{equation}\label{eq:5}
0 = I_0^G \subset I_1^G \subset \cdots \subset I^G_{\dim X}
  \subset I^G_{1 + \dim X} =  A^G
\end{equation}
we only have to show that the algebras in \eqref{eq:10} are in
the category $\ma C$. In fact, since 
$\overline{U_\sigma \setminus D_\sigma} \cap \overline{U_\tau 
\setminus D_\tau} = \emptyset$ if $\dim \sigma = \dim \tau$ and 
$\sigma \neq \tau$, we have an isomorphism
\begin{equation}
I_{m+1} / I_m \cong {\ts \prod_{m + \dim \sigma = \dim X}}
  C_0^\infty ( U_\sigma, D_\sigma ; M_N(\mb C))
\end{equation}
Now $G$ permutes the simplices in this product, so
\begin{equation}\label{eq:11}
I^G_{m+1} / I^G_m \cong {\ts \prod_{\sigma \in L_m}} C_0^\infty 
 (U_\sigma, D_\sigma ; M_N(\mb C))^{G_\sigma}
\end{equation}
where $L_m$ is a set of representatives of the simplices of dimension
$\dim X - m$, modulo the action of $G$. Invoking the additivity of
$K_*$ and $HP_*$ we reduce our task to verifying that every factor
of \eqref{eq:11} belongs to $\ma C$. 
(\emph{This is somewhat problematic! See the appendix.}) 

If $m = \dim X$ then $D_\sigma$ is empty and we see from Lemma 
\ref{lem:7} that $C^\infty (U_\sigma ; M_N(\mb C))^{G_\sigma}$ 
has the required property.

For smaller $m$ there also exists (for every $\sigma$) a 
$G_\sigma$-equivariant contraction of $U_\sigma$ to a point 
$x_\sigma \in \psi(\sigma)$. Thus we can follow the proof of Lemma 
\ref{lem:7} up to equation \eqref{eq:9}, where we find that the 
factor of \eqref{eq:11} corresponding to $\sigma$ is isomorphic to 
$C_0^\infty (U_\sigma, D_\sigma ; \mr{End}_{\mb C}(V_\sigma)
)^{G_\sigma}$. Here $(\pi_\sigma, V_\sigma)$ denotes the projective 
$G_\sigma$-representation over the point $x_\sigma$. Using the 
$G_\sigma$-equivariant retraction of $U_\sigma \setminus D_\sigma$ to 
$\psi(\sigma \setminus D'_\sigma)$ we see that this algebra is 
diffeotopy-equivalent to $C_0^\infty (\sigma, \sigma \cap D'_\sigma) 
\otimes \mr{End}_{G_\sigma} (V_\sigma)$. The right hand side of this
tensor product has finite dimension and is semisimple, so by the 
stability of $\ma C$ it presents no problems. Seen from its 
barycenter $\sigma \setminus D'_\sigma$ is star-shaped, hence by Lemma 
\ref{lem:8} the left hand side is also in the category $\ma C$.

We conclude that all the algebras in \eqref{eq:10} and \eqref{eq:11} 
are indeed objects of $\ma C$, so Lemma \ref{lem:4} can be applied
to \eqref{eq:5} to prove that $A^G \in \mr{Ob}(\ma C)$.

Note that the simplicial complex $\Sigma$ has only finitely many
vertices if $X$ is compact, so then all the above direct 
products are in fact finite and $K_*(A^G)$ is finitely generated.
$\quad \Box$
\\[2mm]

It is clear from the proofs of Lemma \ref{lem:7} and Theorem 
\ref{thm:6} that many similar Fr\'echet algebras are also in $\ma C$. 
For example if $Y$ is a closed submanifold of $X$ then the algebra
\[B = \{f \in C^\infty(X; M_2(\mb C)) : 
  f(y) \: \mr{diagonal} \: \forall y \in Y\}
\]
is in $\ma C$, as can be seen from the admissible extension
\[0 \to C_0^\infty(X,Y;M_2(\mb C)) \to B \to C^\infty(Y)^2 \to 0
\]
The algebra $A$ of Theorem \ref{thm:6} is a finitely generated module
over $C^\infty(X)^G = C^\infty(X/G)$, so more generally one might 
consider Fr\'echet algebras that are finitely generated modules over 
$C^\infty(Y)$ with $Y$ an orbifold. The only real problem to 
generalizing our method to such algebras seems to be that Lemma 
\ref{lem:7} does not apply automatically, so we need a stronger kind
of Poincar\'e lemma. This would require a detailed study of the type
of algebras that can arise in this way.

Now we give some examples of Fr\'echet algebras to which Theorem 
\ref{thm:6} definitely applies.

\begin{cor}\label{cor:9}
Let $\ma H$ be an affine Hecke algebra. Its Schwartz completion 
$\ma S(\ma H)$ belongs to $\ma C$ and has finitely generated
$K$-groups.
\end{cor}
\emph{Proof.}
In Theorem 4.3 of \cite{DO} Delorme and Opdam established an 
isomorphism between $\ma S(\ma H)$ and a finite direct sum of 
algebras of the type described in Theorem \ref{thm:6}, the $X$ 
in each summand being a compact torus. $\quad \Box$
\newpage

Let $G$ be a reductive $p$-adic group. Recall that for a compact 
open subgroup $K$ the Schwartz algebra $\ma S(G//K)$ consists of all 
smooth rapidly decreasing $K$-biinvariant complex valued functions 
on $G$. The convolution product makes it a nuclear Fr\'echet algebra. 
The Schwartz algebra $\ma S(G)$ of $G$ is by definition the union
over all compact open subgroups $\bigcup_K \ma S(G//K)$, 
endowed with the inductive limit topology. It is a complete locally 
convex algebra and a nuclear vector space, but it is not metrizable. 

Furthermore $\ma H(G//K)$ is the subalgebra of $\ma S(G//K)$ 
consisting of compactly supported functions, and $\ma H(G) := 
\bigcup_K \ma H(G//K)$ is called the Hecke algebra of $G$. 
These subalgebras are not complete, and their homologies are 
usually studied with respect to the algebraic tensor product.

Having introduced these objects, we state a crudely simplified 
version of Harish-Chandra's Plancherel formula for $p$-adic 
groups \cite{HC}, a full proof of which was supplied by 
Waldspurger \cite{Wal}.

\begin{thm}\label{thm:10}
There exists a countable collection of triples
$(T_n, L_n, \Gamma_n)$ with the following properties.
For every $n \; T_n$ is a compact torus, $\Gamma_n$ a finite group 
acting on $T_n$ through diffeomorphisms $\alpha_\gamma$ and
$L_n$ an algebra of bounded operators on a Hilbert space. 
There is a $G \times G$-action on $L_n$ such that for any compact 
open subgroup $K$ of $G$ the invariant algebra $L_n^{K \times K}$ 
is semisimple and has finite dimension.
The group $\Gamma_n$ acts on $C^\infty (T_n ; L_n)$ by
\begin{equation}\label{eq:12}
\gamma \cdot f(x) = c_\gamma (x) f(\alpha_\gamma^{-1} x)
\end{equation}
where $c_\gamma (x) \in \mr{Aut}_{G \times G} L_n$. 
All this results in an isomorphism
\begin{equation}\label{eq:13}
\ma S(G) \cong {\ts \bigoplus_{n=1}^\infty} 
  C^\infty (T_n ; L_n)^{\Gamma_n}  
\end{equation}
where the right hand side has the inductive limit topology with 
respect to the direct sum and the $L_n^{K \times K}$.
Furthermore, for every $K$ there exist suitable numbers 
$n_1, \ldots, n_{N_K}$ such that the restriction to 
$K$-biinvariants is an isomorphism
\begin{equation}\label{eq:14}
\ma S(G//K) \cong {\ts \bigoplus_{i=1}^{N_K}} C^\infty (T_{n_i} ; 
  L_{n_i}^{K \times K})^{\Gamma_{n_i}}
\end{equation}
For this action of $\Gamma_{n_i}$ there are  
$u_\gamma \in C^\infty (T_{n_i} ; L_{n_i}^{H \times H})$
such that 
\begin{equation}\label{eq:19}
\gamma \cdot f(x) = u_\gamma (x) f(\alpha_{\gamma}^{-1} x)
u_\gamma^{-1}(x)
\end{equation}
\end{thm}
\emph{Proof.} We will show that deriving this theorem from \cite{Wal}
is merely a matter of translating. We will freely use Waldpurger's
notation, which unfortunately differs substantially from ours, and we 
start by noticing that he writes $\ma C(G)$, respectively
$\ma C_H$, for what we call $\ma S(G)$, respectively $\ma S(G//H)$.
The right hand side of \eqref{eq:13} is 
\[C^\infty (\Theta)^{\mr{inv}} = 
  \big( \bigoplus C^\infty (\ma O, P) \big)^{W^G}
\]
The direct sum runs over all parabolic subgroups $P$ of $G$ that 
contain a certain fixed maximal torus $A_0$. Let $P = M U$ be the
Levi decomposition such that $A_0 \subset M$ and denote the compact 
torus of unitary characters of $M$ by $\mr{Im}\: X(M)$.
Let $(\omega, E)$ be an irreducible square-integrable admissible
representation of $M$ and construct
\[L(\omega, P) := I^{G \times G}_{P \times P} (E \otimes \breve{E}) =
  {\ts \bigcup_{H < G}} L(\omega, P)^{H \times H}
\]
where $I$ denotes induction with respect to compactly supported 
smooth functions, $\breve{E}$ is the contragredient representation
of $E$ and we take the union over all compact open subgroups $H$ of 
$G$. By the admissibility of $E$ each $L (\omega, P)^{H \times H}$ 
has finite dimension, and it is a *-algebra since we can simply 
transfer the * from $\ma S(G//H)$. In particular it is semisimple.

All this leads to the identification
\[C^\infty (\ma O, P) = \big( C^\infty (\mr{Im}\: X(M)) 
  \otimes_{\mb C} L (\omega, P) 
  \big)^{\mr{Stab}_{\mr{Im}\: X(M)} (\omega)}
\]
where the indicated stabilizer acts as in \eqref{eq:12},  
$\alpha_{\gamma}$ being translation by $\gamma$.
Likewise, the Weyl group $W^G$ acts on $C^\infty (\Theta)$ as 
in \eqref{eq:12}, where we should read 
${ }^\circ {c_{P'|P}} (\gamma,x)$ for 
$c_\gamma(x)$. Of course this doesn't fix the $P$'s and the $M$'s, 
so the $c_{\gamma}(x)$ live in something like
$\mr{Hom}_{G \times G} (L_n, L_m)$. In section V.3 it is shown 
that these intertwiners are smooth in $x$.  
Now we take representatives for the association classes of the 
action of $W^G$ on the components of $\Theta$, so that 
$T_n = \mr{Im}\: X(M) \:,\: L_n = L(\omega, P)$ and $\Gamma_n$ 
is semidirect product of the stabilizers of $\omega$ in 
$\mr{Im}\: X(M)$ and of $M$ in $W^G$.

This completes the translation of \eqref{eq:13} to Th\'eor\`emes 
VII.2.5 and VIII.1.1 of \cite{Wal}, leaving \eqref{eq:14} and
\eqref{eq:19}, which are not stated explicitly in that paper.

For an arbitrary compact open subgroup $H$ we consider the 
characteristic function $e_H \in \ma S(G//H)$ defined by
\[e_H (g) = \left\{ \begin{array}{lll}
  0 & \mr{if} & g \notin H \\
  \mr{mes}(H)^{-1} & \mr{if} & g \in H
  \end{array} \right.
\]
It is an idempotent and its image under \eqref{eq:13} lives only in
those components for which
$L(\omega,P)^{H \times H} \neq 0$. These are finite in number and we 
label them by $n_1, \ldots, n_{N_H}$. Since
\[e_H \ma S(G) e_H = \ma S(G//H) = \ma S(G)^{H \times H} 
\]
and the actions of $\Gamma_n$ and $H \times H$ commute, we get 
\eqref{eq:14}. Moreover as $\mr{Aut} (M_N (\mb C)) = PGL (\mb C^N)$,
the automorphism $c_\gamma (x)$ of $L_{n_i}^{H \times H}$ is in fact
conjugation by an invertible element 
$u_\gamma (x) \in L_{n_i}^{H \times H}$ and, $c_\gamma$ being smooth,
we can arrange that 
$u_\gamma \in C^\infty (T_{n_i} ; L_{n_i}^{H \times H}). \qquad \Box$

\begin{cor}\label{cor:11}
$\ma S(G//K) \in \mr{Ob} (\ma C)$ and its $K$-theory is finitely 
generated.
\end{cor}

This corollary provides the small last step needed to complete 
the proof of a conjecture of Baum, Higson and Plymen. Confer 
\cite{BHP}, in particular 8.9 and 9.4, for more background.

\begin{thm}\label{thm:12}
Let $X$ and $C_r^*(G)$ be respectively the affine Bruhat-Tits 
building and the reduced $C^*$-algebra of $G$. 
There exists a commutative diagram
\begin{center}
$\begin{array}{ccc}
K^G_*(X) & \mu \atop \longrightarrow & K_*(C_r^*(G)) \\
\downarrow & & \downarrow \\
HP_*(\ma H(G)) & HP_*(i) \atop \longrightarrow & HP_*(\ma S(G))
\end{array}$
\end{center}
Here the vertical arrows are Chern characters, $\mu$ is the 
Baum-Connes assembly map and $HP_*(i)$ is induced by the inclusion 
$i : \ma H(G) \to \ma S(G)$. The horizontal maps are isomorphisms
and the vertical maps become isomorphisms after tensoring the
diagram with $\mb C$.
\end{thm}
\emph{Proof.}
In \cite{Laf} Lafforgue proved the Baum-Connes conjecture for 
reductive $p$-adic groups, which is another way of saying that $\mu$ 
is an isomorphism. The commutativity of the diagram and the statement 
on the left vertical map can be found in \cite{BHP}. Recall that
$\ma S(G) = \varinjlim \ma S(G//K)$ is a holomorphically closed 
dense subalgebra of $C_r^*(G)$. Hence from the density theorem and 
the continuity of the $K$-functor we get 
\[K_*(C_r^*(G)) \cong \varinjlim K_*(\ma S(G//K))
\] 
To avoid possible problems with the continuity of periodic cyclic 
homology, in \cite{BHP} $HP_*(\ma S(G))$ is defined as 
$\varinjlim HP_*(\ma S(G//K))$. Now Corollary \ref{cor:11}
says that the right vertical map becomes an isomorphism after 
tensoring with $\mb C$, so that the entire diagram will then consist 
of isomorphic objects. Thus $HP_*(i)$, being unmodified by this 
tensoring, is also an isomorphism. $\quad \Box$
\\[2mm]

Actually the formulation of Theorem \ref{thm:12} is somewhat
imprecise since, after specifying a particular topological 
tensor product, we might also calculate $HP_*(\ma S(G))$ directly. 
(Well, in theory at least.) Because $S(G)$ is an inductive limit, it 
is best to use the inductive tensor product. If the result 
would be isomorphic to $\varinjlim HP_*(\ma S(G//K))$, which does not 
seem unlikely, then the theorem has a better and stronger meaning. 

Indeed for $G = GL(m,F)$, with $F$ a non-Archimedean local field,
Theorem \ref{thm:12} was already proved by Brodzki and Plymen 
\cite{BP2}, and they also showed in \cite{BP1} that 
\begin{equation}\label{eq:21}
HP_*(\ma S(G)) \cong \varinjlim HP_*(\ma S(G//K))
\end{equation}
However there are obstacles to generalizing this result to other 
reductive $p$-adic groups. Namely, the proof depends on the vanishing 
of the topological Hochschild homology groups\\ 
$HH_n(\ma S (GL(m,F)//K))$ 
for all $n$ larger then a certain number, independent of $K$. To show
this one uses that $\ma S (GL(m,F)//K)$ is Morita-equivalent to a 
finite direct sum of commutative Fr\'echet algebras of the type 
$C^\infty (X)^W$, for suitable $X$ and $W$. This is a specific 
property of $GL(m,F)$ which does not hold for all other groups.

On the other hand, Nistor \cite{Nis2} showed that
for any reductive $p$-adic group $HH_n(\ma H(G)) = 0$ if $n$ 
exceeds the split rank of $G$. By the continuity of algebraic 
Hochschild homology this means that $HH_n(\ma H(G//K))$ will vanish
for such $n$, for a cofinal collection of compact open subgroups 
$K$. Although Nistor's techniques do not seem to carry over to the 
Schwartz completions, it is not unreasonable to expect that
$HH_n(\ma S(G//K))$ also vanishes for large $n$. 
Using \cite{BP1} this would imply \eqref{eq:21}.

Because of \eqref{eq:25} one might even hope that 
for all the algebras of Theorem \ref{thm:6} we have $HH_n(A^G) = 0$ 
if $n$ exceeds some bound which depends only on $G$ and dim $X$. 
Using the proof of Theorem \ref{thm:6} we can reduce this problem
to the algebras appearing in equation \eqref{eq:11} but there
the applicability of this paper ends, since Lemmas \ref{lem:7} 
and \ref{lem:8} both rely in an essential way on the 
diffeotopy-invariance of $K_*$ and $HP_*$, a property that 
$HH_*$ does not possess.

Finally we remark that recently Meyer \cite{Mey} studied 
the inclusion $\ma H(G) \to \ma S(G)$ from a different perspective, 
namely that of representations of bornological algebras. He obtained
strong results on the comparison of the categories of tempered 
representations of these algebras. It is unclear to the author whether
this has implications for the induced map on periodic cyclic homology.
\newpage

\renewcommand{\theequation}{A.\arabic{equation}}
\renewcommand{\thethm}{A.\arabic{thm}}
\section*{Appendix, added September 2008}

Recently the author became aware that infinite direct products do not commute
with algebraic tensor products. For example, 
$\big( \prod_{n = 1}^\infty \mb Z \big) \otimes_{\mb Z} \mb C$
is strictly smaller than $\prod_{n = 1}^\infty \mb C$. Therefore the map
\[
ch \otimes \mr{id} : K_* \big( {\ts \prod_{n = 1}^\infty \mb C \big) \otimes_{\mb Z} \mb C
\to HP_* \big( \prod_{n = 1}^\infty } \mb C \big)
\]
is injective, but fails to be a surjection. Hence the category $\ma C$
cannot be closed with respect to infinite direct products.

Unfortunately, this very assumption was used in the lines directly following
\eqref{eq:11}. Rather than rewriting a part of the paper, the author decided to add
this appendix, which discusses two ways to overcome the problem.

Firstly, we can avoid it altogether by restricting ourselves to algebras $A$ for which
$HP_* (A)$ has finite dimension. In the context of Theorem \ref{thm:6}, this means that
we have to require that the $G$-manifold $X$ admits a finite open cover, such that 
every intersection of elements of this cover is $G$-equivariantly contractible. This
is the case for compact manifolds, and for algebraic varieties on which 
$G$ acts algebraically.

The second solution is more involved, and replaces $\ma C$ by a somewhat different
class of algebras. The underlying idea is that $\prod_{n = 1}^\infty \mb C$ and
$\big( \prod_{n = 1}^\infty \mb Z \big) \otimes_{\mb Z} \mb C$ are both nondegenerately paired with
$\bigoplus_{n = 1}^\infty \mb C$, and that this property is actually good enough. 
To make sense of this for general $m$-algebras, we use the cohomology theories that 
are dual to $K$-theory and periodic cyclic homology.

Recall that Cuntz and Quillen \cite[Section 10]{CQ1} defined a bivariant functor $HP_* (A,B)$,
such that $HP_* (\mb C ,B) = HP_* (B)$ and $HP_* (A , \mb C ) = HP^* (A)$, where $HP^*$ denotes 
periodic cyclic cohomology. The construction of bivariant periodic cyclic homology can be 
carried out in various categories of algebras, in particular for $m$-algebras. This functor is 
stable and diffeotopy-invariant in both variables \cite[Section 3]{CQ2}. The excision property 
holds with respect to admissible extensions, also in both variables \cite[Section 5]{Cun1}. 
There is a natural product
\[
HP_* (A,B) \times HP_* (B,D) \to HP_* (A,D) \,,
\]
which in particular provides a bilinear form
\begin{equation}\label{eq:A.1}
HP_i (A) \times HP^i (A) \to HP_0 (\mb C ,\mb C ) \cong \mb C \,.
\end{equation}
Moreover Cuntz \cite{Cun2} developed a bivariant $kk$-theory for $m$-algebras, 
whose construction we briefly recall. Let $T A$ be the "projective" completion of the tensor
algebra of $A$, and define $J A$ as the kernel of the multiplication map
$T A \to A$. Let $\mf K$ be the $m$-algebra of infinite matrices with rapidly decaying
coefficients, which is embedded in usual $C^*$-algebra of compact operators. 
Denote the collection of homotopy classes of $m$-algebra homomorphisms $A \to B$ by
$\inp{A}{B}$. There exist natural maps $\inp{J^n A}{B} \to \inp{J^{n+2} A}{\mf K \hot B}$, 
where $\hot$ denotes the completed projective tensor product. Cuntz defines
\begin{equation}\label{eq:A.8}
kk_j (A,B) = \lim_{n \to \infty} \inp{J^{2n+j} A}{\mf K \hot B} \,.
\end{equation}
This theory takes values in abelian groups and is 2-periodic. Moreover $kk_*$ satisfies all the
formal properties of $HP_*$ which we described above. We write
\[
K_* (A) = kk_* (\mb C ,A) \quad \mr{and} \quad K^* (A) = kk_* (A,\mb C) \,. 
\]
This agrees with the classical $K$-theory for Banach algebras, and with Phillips' \cite{Phi}
$K$-theory for Fr\'echet algebras. However, it differs from Kasparov's $KK$-theory and from 
$K$-homology for separable $C^*$-algebras.
The product in $kk$-theory gives a $\mb Z$-bilinear form
\begin{equation}\label{eq:A.2}
K_i (A) \times K^i (A) \to kk_0 (\mb C ,\mb C ) \cong \mb Z \,.
\end{equation}
Finally, there exists a functorial bivariant Chern character
\[
ch : kk_* (A,B) \to HP_* (A,B) \,,
\]
which is compatible with all these properties.

In general $kk_* (\prod_{n = 1}^\infty A_n , \prod_{m = 1}^\infty B_m )$ is neither isomorphic to
$\prod_{m = 1}^\infty kk_* (\prod_{n = 1}^\infty A_n , B_m )$, nor to 
$\bigoplus_{n = 1}^\infty kk_* (A_n , \prod_{m = 1}^\infty B_m )$.
For example, the identity map $\mr{id}_P$ of $P = \prod_{n = 1}^\infty \mb C$ 
cannot be represented by an element of $\bigoplus_{n = 1}^\infty kk_* (\mb C , P)$.
Similar phenomena occur in bivariant periodic cyclic homology.
Fortunately the univariant versions of these functors do behave well 
with respect to infinite direct products:

\begin{prop}\label{prop:A.1}
Let $A_i \; (i \in I)$ be an arbitrary collection of $m$-algebras.
There are natural isomorphisms
\begin{align}
\label{eq:A.3} HP_* \big( {\ts \prod_{i \in I}} A_i \big) & \cong {\ts \prod_{i \in I} } HP_* (A_i ) \,, \\
\label{eq:A.4} HP^* \big( {\ts \prod_{i \in I}} A_i \big) & \cong {\ts \bigoplus_{i \in I} } HP^* (A_i ) \,, \\
\label{eq:A.5} K_* \big( {\ts \prod_{i \in I}} A_i \big) & \cong {\ts \prod_{i \in I} } K_* (A_i ) \,, \\
\label{eq:A.6} K^* \big( {\ts \prod_{i \in I}} A_i \big) & \cong {\ts \bigoplus_{i \in I} } K^* (A_i ) \,.
\end{align}
\end{prop}
\emph{Proof.}
Let $B$ be any $m$-algebra and abbreviate $A = \prod_{i \in I} A_i$.\\
By definition $HP_* (B)$ is the homology of a differential complex $\ma X (\overleftarrow{T} B)$, see
\cite[Section 4.1.5]{Mey2}. According to \cite[Theorem 4.3.7]{Mey2} the canonical map
\begin{equation}\label{eq:A.7}
\ma X (\overleftarrow{T} A) \to {\ts \prod_{i \in I}} \ma X (\overleftarrow{T} A_i )
\end{equation}
is a homotopy equivalence, which leads to \eqref{eq:A.3}.

Furthermore $HP^* (B)$ is the cohomology of Hom$(\ma X (\overleftarrow{T} B), \mb C)$, where we must
take the homomorphisms in the category of projective limits of Banach spaces. Recall that the continuous
linear dual of a direct product of topological vector spaces can be identified with the direct sum of
the dual spaces. Together with \eqref{eq:A.7} we find homotopy equivalences
\[
\mr{Hom}(\ma X (\overleftarrow{T} A), \mb C) \leftarrow 
\mr{Hom} \big( {\ts \prod_{i \in I} \ma X (\overleftarrow{T} A_i ), \mb C \big) \cong 
\bigoplus_{i \in I} } \mr{Hom}( \ma X (\overleftarrow{T} A_i ), \mb C ) \,,
\]
which implies \eqref{eq:A.4}.

From \cite[p. 178]{Cun2} we see that
\[
\begin{array}{ccccl}
K_0 (B) & = & kk_0 (\mb C ,B) & \cong & \inp{q \mb C}{\mf K \hot B} \\
K_1 (B) & = & kk_1 (\mb C ,B) & \cong & \inp{q \mb C}{\mf K \hot \mb C (0,1) \hot B}  \,,
\end{array}
\]
where the Fr\'echet algebras $q \mb C$ and $\mb C (0,1)$ are as in \cite[Section 1]{Cun2}.
Using the compatibility of $\hot$ with direct products \cite[Proposition I.1.3.6]{Gro1} we get
natural isomorphisms
\[
{\ts K_0 \big( \prod_{i \in I} A_i \big) \cong \inp{q \mb C}{\prod_{i \in I} (\mf K \hot A_i )} 
\cong \prod_{i \in I} \inp{q \mb C}{\mf K \hot A_i} \cong \prod_{i \in I} K_0 (A_i ) } \,.
\]
The same goes for $K_1 (A)$, proving \eqref{eq:A.5}.

On the other hand, for $B = \mb C$ \eqref{eq:A.8} becomes
\[
K^j (A) = kk_j (A,\mb C) = \lim_{n \to \infty} \inp{J^{2n+j} A}{\mf K} \,.
\]
For any finite subset $F \subset I$ we write $A_F = \bigoplus_{i \in F} A_i$. 
Since $kk_*$ satisfies excision, 
\[
K^* (A_F ) \cong {\ts \bigoplus_{i \in F}} K^* (A_i ) \,.
\]
The inclusion and quotient maps $A_F \to A \to A_F$ induce group homomorphisms
\[
K^* (A_F ) \to K^* (A) \to K^* (A_F ) \,,
\]
whose composition is the identity. These combine to a natural injection
\begin{equation}\label{eq:A.9}
{\ts \bigoplus_{i \in I}} K^* (A_i ) \to K^* (A) \,.
\end{equation}
Consider the following subalgebra of $J^m (A)$:
\[
(J^m A)_F := J^m (A_F ) \cap {\ts \bigcap_{i \in F}} 
\ker \big( J^m (A_F ) \to J^m \big( A_{F \setminus \{ i \}} \big) \big) \,.
\]
It can be described as the subspace of $J^m (A)$ which is spanned by all tensors which only 
involve elements from the $A_i$ with $i \in F$, and which is complementary to the tensors
coming from fewer summands $A_i$. Notice that
\[
J^m (A_F ) = {\ts \prod_{F' \subset F}} (J^m A )_{F'} \,.
\]  
This product is direct in the category of topological vector spaces, but as a product of 
algebras it is only semi-direct.
Since $J^m A = \varprojlim J^m (A_F )$, we can identify it as a topological vector space with
the direct product of the $(J^m A )_F$, over all finite subsets $F$ of $I$:
\begin{equation}\label{eq:A.11}
J^m (A) = {\ts \prod_F} (J^m A )_F \,.
\end{equation}
To show that \eqref{eq:A.9} is surjective, take any class $[f] \in K^j (A)$, 
represented by an $m$-algebra homomorphism $f : J^{2n + j} A \to \mf K$. 
Let $\| ? \|_o$ be operator norm on the pre-$C^*$-algebra $\mf K$, so that
\[
J^{2n+j} A \to \mb R \,:\, a \mapsto \| f(a) \|_o
\]
is a continuous map. Thus
\[
U := \{ a \in J^{2n+j} A : \| f(a) \|_o < 1 \}
\] 
is open. Since the right hand side of \eqref{eq:A.11} is equipped with the product topology, 
we have $(J^{2n+j} A )_F \subset U$ for all but finitely many $F$.
Clearly this implies $f \big( (J^{2n+j} A )_F \big) = 0$ for these $F \subset I$. Let $I_f$
be the union of the $F \subset I$ for which $f \big( (J^{2n+j} A )_F \big) \neq 0$. Then
$I_f$ is finite and $f$ factors as
\[
J^{2n+j}(A) \to J^{2n+j} \big( A_{I_f} \big) \to \mf K \,.
\]
This means that $[f]$ lies in the image of $K^j \big( A_{I_f} \big) \to K^j (A)$, so
\eqref{eq:A.9} is indeed surjective. $\qquad \Box$
\\[2mm]

In many cases the pairings \eqref{eq:A.1} and \eqref{eq:A.2} are nondegenerate, so it
makes sense to consider the following class of algebras.

\begin{defn}\label{def:A.2}
The class $\ma C'$ consists of all $m$-algebras $A$ satisfying the following conditions:
\begin{enumerate}
\item $HP_* (A) \times HP^* (A) \to \mb C$ is a nondegenerate bilinear pairing,
\item $ch : K^* (A) \otimes_{\mb Z} \mb C \to HP^* (A)$ is a linear bijection,
\item $K_* (A) \otimes_{\mb Z} \mb C \times K^* (A) \otimes_{\mb Z} \mb C \to  \mb C$ 
is a nondegenerate bilinear pairing.
\end{enumerate}
\end{defn}

\noindent For any $m$-algebra $A$ in $\ma C'$, the Chern character 
\begin{equation}\label{eq:A.10}
ch : K_* (A) \otimes_{\mb Z} \mb C \to HP_* (A)
\end{equation}
is injective and has dense image, with respect to the coarsest topology on $HP_* (A)$
that makes all elements of $HP^* (A)$ into continuous linear functionals.
In particular, if $HP_* (A)$ has finite dimension, then \eqref{eq:A.10} is a linear 
bijection and $A$ also belongs to $\ma C$.

From the above discussion we already know that $\ma C'$ is closed under diffeotopy
equivalences and under tensoring with $M_k (\mb C )$. Proposition \ref{prop:A.1} and
the compatibility of $\otimes_{\mb Z} \mb C$ with direct sums show that $\ma C'$
is closed with respect to arbitrarily large direct products.

After these remarks we will
show that $\ma C'$ has the same "two out of three"-property as $\ma C$.
Hence $\ma C'$ really has all the properties claimed for $\ma C$ in Corollary \ref{cor:3}. 
Since these are all that is used in the remainder of the paper, we conclude that all the 
results, in particular Theorem \ref{thm:6}, hold if we replace $\ma C$ by $\ma C'$. 
Moreover Corollaries \ref{cor:9} and \ref{cor:11} are valid precisely as stated (with $\ma C$), 
because the periodic cyclic homology has finite dimension in these cases.

\begin{lem}\label{lem:A.3}
Let $0 \to A \xrightarrow{\phi} B \xrightarrow{\psi} D \to 0$ be an admissible extension of 
$m$-algebras. If two elements of $\{ A,B,D \}$ belong to $\ma C'$, then so does the third. 
\end{lem}
\emph{Proof.}
Of the three conditions in Definition \ref{def:A.2}, we can handle the second just as in
the proof of Corollary \ref{cor:3}. As moreover all cases of the first and third conditions
can be dealt with in the same way, we will only treat one case. Suppose that $A$ and $D$
belong to $\ma C'$. We want to show that $HP_0 (B) \times HP^0 (B) \to \mb C$ is a 
nondegenerate bilinear pairing. Consider the exact sequences
\[
\begin{array}{*{13}{c}}
\cdots & \to & HP_1 (D) & \to & HP_0 (A) & \to & HP_0 (B) & \to & 
HP_0 (D) & \to & HP_1 (A) & \to & \cdots \\
\cdots & \leftarrow & HP^1 (D) & \leftarrow & HP^0 (A) & \leftarrow & HP^0 (B) &
\leftarrow & HP^0 (D) & \leftarrow & HP^1 (A) & \leftarrow & \cdots
\end{array}
\]
In every column there is a bilinear pairing, and every arrow is adjoint to the one
directly below or above it. By assumption all these pairings are nondegenerate, except
possibly in the middle column. Therefore
\[
\begin{array}{l@{\;=}l@{\qquad}l@{\;=}l}
\mr{im}(HP_1 (D) )^\perp & \ \; \mr{im} (HP^0 (B)) & 
\mr{im} (HP_0 (B) )^\perp & \ \; \mr{im}(HP^1 (A)) \,, \\
\mr{im}(HP_1 (D) )& \ \; \mr{im} (HP^0 (B) )^\perp & 
\mr{im} (HP_0 (B) ) & \ \; \mr{im}(HP^1 (A) )^\perp ,
\end{array}
\]
and we can simplify the above diagram to
\[
\begin{array}{ccccccccc}
0 & \to & HP_0 (A) / \mr{im}(HP_1 (D)) & \to & HP_0 (B) & \to & \mr{im}(HP_0 (B)) & \to & 0 \\
0 & \leftarrow & \mr{im}(HP^0 (B)) & \leftarrow & HP^0 (B) & \leftarrow & 
HP^0 (D) / \mr{im}(HP^1 (A)) & \leftarrow & 0
\end{array}
\]
The rows remain exact and the second and fourth columns are endowed with nondegenerate 
bilinear pairings. Consider any $x \in HP_0 (B) \cap HP^0 (B)^\perp$. Then 
$HP_0 (\psi )(x) \in HP^0 (D)^\perp$, so by the nondegeneracy of 
\[
\mr{im} (HP_0 (B)) \times HP^0 (D) / \mr{im}(HP^1 (A)) \to \mb C
\]
we have $HP_0 (\psi )(x) = 0$. Thus $x = HP_0 (\phi )(y)$ for some $y \in HP_0 (A)$. But 
\[
y \in \mr{im}(HP^0 (B) )^\perp = \mr{im}(HP_1 (D)) , \text{ so } x \in HP_0 (\phi) (HP_1 (D)) = 0 .
\]
A similar argument shows that $HP^0 (B) \cap HP_0 (B)^\perp = 0. \qquad \Box$

\vspace{1cm}

\noindent \textbf{Acknowledgements.}

The author would like to thank Ralf Meyer, Victor Nistor, Eric Opdam and Peter 
Schneider for their helpful comments.

\end{document}